\documentclass[aap,MSNbibl,seceqn,dvips]{arximspdf}

\doi{10.1214/13-AAP922} 
\volume{24}
\issue{2}
\pubyear{2014}
\firstpage{526}
\lastpage{552}

\makeatletter
\newcommand{\rrvert}{\vert}
\newcommand{\llvert}{\vert}
\newcommand{\II}{\mathrm{I}}
\newcommand{\Tr}{\operatorname{Tr}}
\newcommand{\Law}{\operatorname{Law}}
\def\E{\mathrm{E}}
\def\P{\mathrm{P}}
\newtheorem{Theorem} {Theorem}[section]
\newtheorem{Lemma}[Theorem] {Lemma}

\newproclaim{Example}[Theorem]{Example}
\newproclaim{Remark}[Theorem]{Remark}
\newproclaim{Definition}{Definition}[section]
\makeatother

\begin{document}
\begin{frontmatter}

\title{Subgeometric rates of convergence of Markov processes in the
Wasserstein metric}
\runtitle{Rates of convergence in the Wasserstein metric}

\begin{aug}
\author[A]{\fnms{Oleg} \snm{Butkovsky}\corref{}\ead[label=e1]{oleg.butkovskiy@gmail.com}\thanksref{t1}}
\thankstext{t1}{Supported in part by Russian Foundation for Basic Research
Grant 10-01-00397-a, Israel Science Foundation Grant 497/10 and a Technion fellowship.}
\runauthor{O. Butkovsky}
\affiliation{Lomonosov Moscow State University and\\ Technion---Israel Institute of Technology}
\address[A]{Faculty of Mathematics and Mechanics\\
Department of Probability Theory\\
Lomonosov Moscow State University\\
Moscow, 119991\\
Russia\\
and\\
Faculty of Industrial Engineering and Management\\
Technion---Israel Institute of Technology\\
Haifa 32000\\
Israel\\
\printead{e1}} 
\end{aug}

\received{\smonth{11} \syear{2012}}
\revised{\smonth{1} \syear{2013}}

%
\begin{abstract}
We establish subgeometric bounds on convergence rate of general Markov
processes in the Wasserstein metric. In the discrete time setting we
prove that the Lyapunov drift condition and the existence of a ``good''
$d$-small set imply subgeometric convergence to the invariant measure.
In the continuous time setting we obtain the same convergence rate
provided that there exists a ``good'' $d$-small set and the
Douc--Fort--Guillin supermartingale condition holds. As an application
of our results, we prove that the Veretennikov--Khasminskii condition
is sufficient for subexponential convergence of strong solutions of
stochastic delay differential equations.
\end{abstract}

%

\end{frontmatter}

\section{Introduction}\label{sec1}

In this paper, we study rate of convergence of Markov processes to an
invariant measure in the Wasserstein metric. We establish subgeometric
bounds on the convergence rate, thus generalizing the results of
\cite{DFG,DFMS,HMS}. We apply the obtained estimates to prove
subgeometric ergodicity of strong solutions of stochastic differential
delay equations (SDDEs) under Veretennikov--Khasminskii-type
conditions. This extends the corresponding results
\cite{Ver97,VerTVP,Mal,DFG} for stochastic differential equations
(without delay).

There are quite a few works which deal with convergence of Harris
recurrent Markov chains in total variation; see, for example, the
monograph \cite{MT} and the references therein. Less is known about
convergence of Markov chains that are not Harris recurrent. Recall
\cite{HLL} that if a Markov chain has a unique invariant measure, then
either (a) the chain is positive Harris recurrent in an absorbing set
and the invariant measure is nonsingular, or (b) the invariant measure
is singular and there are no Harris sets. It is quite clear that in
case~(b) the marginal distributions of the Markov chain do not converge
in total variation, whereas they might converge weakly (and, hence, in
the Wasserstein metric). Thus, for non-Harris chains [case~(b)] it is
natural to study convergence in the Wasserstein metric (rather than in
the total variation metric).

Many interesting Markov processes fall into case (b). For instance,
following~\cite{HMS}, consider SDDE
\[
d X(t) = - c X(t)\,dt+g \bigl(X(t-1) \bigr)\,dW(t),\qquad t>0,
\]
where $c>0$, $W$ is a one-dimensional Brownian motion and $g$ is a
strictly increasing positive bounded continuous function. One can show
that the strong solution of this equation has a unique invariant
measure and converges to it weakly, but not in total variation. On the
other hand, the Wasserstein distance between $X(t)$ and the invariant
measure decays exponentially to zero as $t\to\infty$. Section~\ref{S3}
contains further examples of processes belonging to case~(b).

Many methods of estimation of convergence rates in the total variation
metric assume that a Markov process is $\psi$-irreducible and are based
on the analysis of small sets. Probably, one of the first results in
this area is due to Dobrushin \cite{Dobr}, who proved that if the whole
state space is small, then a Markov chain is exponentially ergodic.
Later Popov \cite{Pop} and Nummelin and Tuominen \cite{NT82} replaced
the global Dobrushin condition with a combination of a local Dobrushin
condition (existence of a ``good'' small set) and the Lyapunov drift
condition (LDC). This result was further extended by Jarner and Roberts
\cite{JR} and Douc and coauthors \cite{DFMS}, who established
polynomial and general subgeometric estimates of convergence rate,
correspondingly. Similar results for continuous time Markov processes
(under an additional assumption that the state space is locally
compact) are due to Fort and Roberts \cite{FR} and Douc, Fort and
Guillin \cite{DFG}. The latter work provides subgeometric estimates of
the convergence rate under condition that a certain functional of a
Markov process is a supermartingale. Let us also mention the recent
paper of Hairer and Mattingly \cite{HM}, which contains a new simple
proof of the exponential ergodicity of a Markov process under LDC and
the local Dobrushin condition.

Thus, many techniques rely on the irreducibility of a Markov process,
the existence of a ``good'' small set, and (for continuous time
processes) the local compactness of the state space. However, if the
state space is infinite-dimensional, then in most ``typical''
situations the process is non-Harris and, therefore these assumptions
are not fulfilled. For instance, if we go back to the above SDDE, then
it is easy to check that this processes is not $\psi$-irreducible, the
state space is not locally compact and, as was pointed in \cite{HMS},
all small sets of this process are degenerate (i.e., consists of no
more than one point).

An alternative to the local Dobrushin condition was suggested by Bakry,
Cattiaux and Guillin in \cite{BCG}. They obtained estimates of
convergence rate in the total variation metric, provided that the LDC
holds, and a Markov process has a unique invariant measure, which
satisfies a local Poincar\'{e} inequality on a large enough set.

Let us discuss another alternative to this set of assumptions, which
was developed by Hairer, Mattingly, and Scheutzow \cite{HMS}
specifically for establishing exponential convergence rates of SDDEs,
stochastic PDEs, and other infinite-dimensional processes in the
Wasserstein metric. Exploiting a new notion of a $d$-small set
(a~generalization of the notion of a small set), in conjunction with
the LDC, and without any additional assumptions on the irreducibility
of the process, the authors proved the existence of a spectral gap in a
suitable norm, and, hence, the exponential convergence to
stationarity.\looseness=-1

We extend this result and consider the more general situation where a
spectral gap may not exist. For discrete time Markov processes (Theorem~\ref{ThDiscreteTime})
we prove that existence of a ``good'' $d$-small
set and the LDC implies subgeometrical convergence in the Wasserstein
metric. In the continuous time setting (Theorem~\ref{ThContinuousTime})
we obtain the same rate of convergence provided that there exists a
``good'' \mbox{$d$-}small set and the Douc--Fort--Guillin
supermartingale condition holds. Thus, we also extend the results of
\cite{DFG,DFMS}.

We apply our conditions to study the asymptotic behavior of strong
solutions of SDDEs. We prove that Veretennikov--Khasminskii-type
conditions are sufficient for subexponential ergodicity (Theorem
\ref{ThSDDEcoeff}). This extends the results
of~\mbox{\cite{Ver97,VerTVP,Mal,DFG}}.

The rest of the paper is organized as follows. Section~\ref{S2}
contains definitions and the main results. Applications to SDDEs and to
an autoregressive model are presented in Section~\ref{S3}. The proofs
of the main results are placed in Section~\ref{sect4}.

\section{Main results}\label{S2}

Let $X=(X_n)_{n\in\mathbb{Z}_+}$ be a homogeneous Markov chain on a
measurable space $(E,\mathcal{B}(E))$ with transition functions
$P^n(x,A):=P_x(X_n\in A)$, where $x\in E$, $A\in\mathcal{B}(E)$,
$n\in\mathbb{Z}_+$. As usual for $n=1$ we will drop the upper index and
write $P(x,A)$. For a measurable function $f\dvtx E\to[0,\infty)$, let
$\mathcal{P}_f(E)$ be the set of probability measures on
$(E,\mathcal{B}(E))$ which integrate $f$. We will write
$\mathcal{P}(E)$ for the set of all probability measures on
$(E,\mathcal{B}(E))$. If $\mu\in\mathcal{P}_f(E)$, denote
$\mu(f):=\int_E f(x) \mu(dx)$. We define Markov semigroup operators as
usual,
\[
P\varphi(x):=\int_E \varphi(t) P(x,dt),\qquad P \mu(dx):=
\int_E P(t,dx) \mu(dt).
\]

Recall (see, e.g., \cite{BK}) that if $d$ is a semimetric on $E$, then
the \textit{Wasserstein semidistance} $W_d$ between probability
measures $\mu,\nu\in\mathcal{P}(E)$ is given by
\[
W_d(\mu,\nu):=\inf_{\lambda\in\mathcal{C}(\mu,\nu)} \int_{E\times E}d(x,y) \lambda(dx,dy),
\]
where $\mathcal{C}(\mu,\nu)$ is the set of all probability measures
on $(E\times E,\mathcal{B}(E\times E))$ with
marginals $\mu$ and $\nu$. If $d$ is a proper metric, then $W_d$ is a
distance.

We consider also the \textit{total variation metric} on the space
$\mathcal{P}(E)$, which is defined by the following formula:
\[
d_{\mathrm{TV}}(\mu,\nu):=2\sup_{A\in\mathcal{B}(E)} \bigl|\mu(A)-\nu(A)\bigr|,
\qquad\mu,\nu\in\mathcal{P}(E).
\]
Recall that if the space $E$ is equipped with the discrete metric
$d_0(x,y):=\II(x\neq y)$, $x,y\in E$, then the Wasserstein distance is
just half of the total variation distance, that is,
$W_{d_0}(\mu,\nu)=d_{\mathrm{TV}}(\mu,\nu)/2$,
$\mu,\nu\in\mathcal{P}(E)$.

%
\begin{Definition}
A set $A\in\mathcal{B}(E)$ is called \textit{small} for a Markov
operator $P$ if there exists $\varepsilon>0$ such that for all $x,y\in
A$,
\[
\tfrac12d_{\mathrm{TV}} \bigl(P(x,\cdot),P(y,\cdot) \bigr)\leq1-\varepsilon.
\]
\end{Definition}

For instance, any one-point set is small. However, as discussed above,
a Markov process might have no small sets that consist of more than one
point. To study such Markov processes Hairer, Mattingly and Scheutzow
\cite{HMS} introduce the following concept.

%
\begin{Definition}
A set $A\in\mathcal{B}(E)$ is called \textit{$d$-small} for a Markov
operator $P$ if there exists $\varepsilon>0$ such that for all $x,y\in
A$,
\[
W_d \bigl(P(x,\cdot),P(y,\cdot) \bigr)\leq(1-\varepsilon)\,d(x,y).
\]
\end{Definition}

Note that our definition of a $d$-small set is a bit different from the
definition of~\cite{HMS}. Namely, the multiplier $d(x,y)$ appears on
the right-hand side of the above inequality.

If $d(x,y) = \II(x\neq y)$, then the notions of a small set and a
$d$-small set coincide. In the general case, the latter notion is much
weaker than the former. In Section~\ref{Sect31} we give an example of
a Markov operator $P$ that has a $d$-small state space and no
nontrivial small sets.

Before we present our main result, let us recall that the total
variation metric is contracting, that is, for any
Markov semigroup $(P^t)_{t\geq0}$ one has
\[
d_{\mathrm{TV}} \bigl(P^t(x,\cdot),P^t(y,\cdot) \bigr)
\leq d_{\mathrm{TV}} \bigl(P^s(x,\cdot),P^s(y,\cdot)
\bigr),\qquad x,y\in E
\]
whenever $0\leq s\leq t$. In general, the Wasserstein metric $W_d$ may
not be contracting. However, as discussed in detail in \cite{HMS}, it
is natural to focus only on Wasserstein metrics that are contracting
for the process $X$, since, in the general case, the Lyapunov drift
condition is not sufficient even for a weak convergence toward the
invariant measure. Note that the contractivity condition itself does
not imply any convergence at all, either. It is the combination of the
contractivity, the Lyapunov drift condition and the existence of a
``good'' $d$-small set, which yields the existence and uniqueness of
the invariant measure and subgeometric convergence in the Wasserstein
metric.

For a function $f\dvtx\mathbb{R}_+\to(0;\infty)$ define
\[
H_f(x):=\int_1^x\frac{1}{f(u)}
\,du,\qquad x\geq1.
\]
Since $H_f$ is increasing, the inverse function $H_f^{-1}$ is well defined.

\begin{Theorem}\label{ThDiscreteTime}
Suppose there exist a measurable function $V\dvtx E\to[0;\infty)$ and a
metric $d$ on $E$ such that the following conditions hold:
\begin{longlist}[(3)]
\item[(1)] $V$ is a Lyapunov function; that is, there exist a
concave differentiable function $\varphi
\dvtx\mathbb{R}_+\to\mathbb{R}_+$ increasing to infinity with
$\varphi(0)=0$ and a constant $K\geq0$ such that
%
%
\begin{equation}
\label{MTLyapunov} PV\leq V- \varphi\circ V + K.
\end{equation}

\item[(2)] The space $(E,d)$ is a complete separable metric space.

\item[(3)] The metric $d$ is contracting and bounded by $1$; that is,
for any $x,y\in E$,
%
%
\begin{equation}
\label{MTdcontr} W_d \bigl(P(x,\cdot),P(y,\cdot) \bigr)\leq d(x,y)
\leq1.
\end{equation}

\item[(4)] The level set $L:=\{x,y\in E\dvtx V(x)+V(y)\leq R\}$ is
$d$-small for some $R>\varphi^{-1}(2K)$; that is, there exists
$\rho>0$ such that
\[
W_d \bigl(P(x,\cdot),P(y,\cdot) \bigr)\leq(1-\rho)\,d(x,y)
\]
for any $x,y\in L$.
\end{longlist}

Then the process $X$ has a unique stationary measure $\pi$ and
\[
\int_E \varphi \bigl(V(u) \bigr) \pi(du)\leq K.
\]
Moreover, for any $\varepsilon>0$ there exist constants $C_1$ and
$C_2$ such that for all $x\in E$,
%
%
\begin{equation}
\label{convrate} W_d \bigl(P^{n}(x,\cdot),\pi \bigr)\leq
\frac{C_1 (1+V(x))}{ \varphi(H_\varphi
^{-1}(C_2n))^{1-\varepsilon}}, \qquad n\in\mathbb{Z}_+.
\end{equation}
\end{Theorem}

%
\begin{Remark}
(i) If $\varphi$ is a linear function, then the rate of convergence is
exponential and this case is covered by \cite{HMS}, Theorem 4.8.

(ii) If $d(x,y)=\II(x\neq y)$, then the Wasserstein metric coincides
with the total variation metric and this case is covered by
\cite{DFMS}, Proposition 2.5.
\end{Remark}

%
\begin{Remark}
Conditions (3) and (4) of the theorem are a bit more general than the
corresponding conditions from \cite{HMS}, Theorem 4.8. Namely, we do
not assume here that $W_d(P(x,\cdot),P(y,\cdot))\leq(1-\rho)\,d(x,y) $
for all $x,y\in E$ such that $d(x,y)<1$. We suppose that this
inequality is satisfied only for $x$, $y$ belonging to the sublevel
set.
\end{Remark}

Note that if $\varphi$ grows to infinity not very rapidly (as
$x^\gamma$ for some $0<\gamma<1$ or slower), then the estimate of
convergence rate given by (\ref{convrate}) can be as close as possible
to the estimate of convergence rate in the total variation distance
obtained in \cite{DFMS}, Proposition 2.5. Specific examples of
convergence rates (polynomial, logarithmic, etc.) for different
functions $\varphi$ are given in \cite{DFMS}, Section~2.3.

While the proof of the theorem is postponed to Section~\ref{sect4}, we
outline now the main steps.

\begin{pf*}{Sketch of the proof of Theorem~\ref{ThDiscreteTime}}
To prove the theorem we develop the idea of constructing an auxiliary
contracting semimetric \cite{Hair,HM,HMS}. Namely, let $l$ be a
semimetric on the space $E$ such that $d(x,y)\leq l(x,y)$ for all
$x,y\in E$. It is possible to prove (for some ``good'' $l$) that for
any probability measures $\mu,\nu\in\mathcal{P}_{\varphi\circ V}(E)$
\[
W_l(P\mu,P\nu)\leq \bigl(1-\chi(\mu,\nu) \bigr)W_l(
\mu,\nu),
\]
where $\chi$ is a positive function (this is done in Lemma~\ref{LL3}).
Hence
\[
W_d \bigl(P^n\mu,P^n\nu \bigr)\leq
W_l \bigl(P^n\mu,P^n\nu \bigr)\leq\prod
_{i=0}^{n-1} \bigl(1-\chi
\bigl(P^i \mu,P^i \nu \bigr) \bigr)W_l(\mu,
\nu).
\]
Of course, since we want to obtain subgeometric estimates of
$W_d(P^n\mu,P^n\nu)$, there is no hope that
$\inf_{\mu,\nu\in\mathcal{P}_{\varphi\circ V}(E)} \chi(\mu,\nu)$ is
positive (this lower bound was greater than zero in \cite{Hair,HM,HMS},
where geometric estimates were obtained). Yet, a~good (albeit
nonuniform) estimate of $\chi(P^{i+1}\mu,P^{i+1} \nu)$ can be derived.
However, this estimate depends not only on $W_l(P^{i}\mu,P^{i}\nu)$ but
also on $\mu(P^i(\varphi\circ V))$ and $\nu(P^i(\varphi\circ V))$.
The latter two expressions are unbounded if $\mu,\nu$ are fixed, and
$i$ runs over positive integers. Fortunately, there are sufficiently
many integers $i$ such that these two expressions are ``small'' (Lemma
\ref{LL1}). This allows us to overcome this obstacle (Lemma~\ref{LL4})
and obtain subgeometric bounds on $W_d(P^n\mu,P^n\nu)$. The last step
is to prove the existence and uniqueness of the stationary measure
(Lemma~\ref{LL5}).
\end{pf*}

Now we give a similar result for continuous time Markov processes. Let
$X=(X_t)_{t\geq0}$ be a time-homogeneous strong Markov process, and let
$(P_t)_{t\geq0}$ be the associated Markov semigroup. Recall \cite{DY},
Theorem 2, that if a Markov process has c\`{a}dl\`{a}g paths, then the
strong Markov property is implied by the Feller property.

\begin{Theorem}\label{ThContinuousTime}
Suppose there exist a measurable function $V\dvtx E\to[0;\infty)$ and a
metric $d$ on $E$ such that the following conditions hold:
\begin{longlist}[(2)]
\item[(1)] $V$ is a Lyapunov function; that is, there exist a
concave differentiable function $\varphi\dvtx
\mathbb{R}_+\to\mathbb{R}_+$ increasing to infinity with\vadjust{\goodbreak}
$\varphi(0)=0$ and a constant $K\geq0$ such that for all
$t\geq0$, $x\in E$
%
%
\begin{equation}
\label{CTintegr} \E_x V(X_t)\leq V(x) - \E_x
\int_0^t \varphi \bigl(V(X_u)
\bigr)\,du + Kt.
\end{equation}

\item[(2)] The space $(E,d)$ is a complete separable metric space.

\item[(3)] The metric $d$ is bounded by $1$ and contracting for all
$t\geq t_0$, for some $t_0\geq0$; that is, for any $x,y\in E$
\[
W_d \bigl(P^{t}(x,\cdot),P^{t}(y,\cdot)
\bigr) \leq d(x,y)\leq1.
\]

\item[(4)] The level set $L:=\{x,y\in E\dvtx V(x)+V(y)\leq R\}$ is
$d$-small for all $R>0$ and all $t\geq t_0$, that is, there
exists $\rho=\rho(R,t)>0$ such that
\[
W_d \bigl(P^t(x,\cdot),P^t(y,\cdot) \bigr)
\leq(1-\rho)\,d(x,y)
\]
for any $x,y\in L$.
\end{longlist}

Then the process $X$ has a unique stationary measure $\pi$ and $\pi
(\varphi\circ V)\leq K$. Moreover, for any
$\varepsilon>0$ there exist constants $C_1$ and $C_2$ such that for
all $x\in E$,
%
%
\begin{equation}
\label{CTrate} W_d \bigl(P^{t}(x,\cdot),\pi \bigr)\leq
\frac{C_1 (1+V(x))}{ \varphi(H_\varphi
^{-1}(C_2t))^{1-\varepsilon}}, \qquad t\geq0.
\end{equation}
\end{Theorem}

%
\begin{Remark}
(i) The linear case $\varphi(x)=\lambda x$, $\lambda
>0$ is \cite{HMS}, Theorem 4.8.

(ii) The case where the metric $d$ is discrete, that is,
$d(x,y)=\II(x\neq y)$, is \cite{DFG}, Theorem 3.2.
\end{Remark}

%
\begin{Remark}
(i) Condition (1) of Theorem~\ref{ThContinuousTime} is equivalent to the
Douc--Fort--Guillin supermartingale condition \cite{DFG},
equation~(3.2); that is, inequality (\ref{CTintegr}) holds if and only
if the process $Z:=(Z_t)_{t\geq0}$,
\[
Z_t:=V(X_t)+\int_0^t
\varphi \bigl(V(X_u) \bigr)\,du-Kt,\qquad t\geq0
\]
is a supermartingale with respect to the natural filtration of the
process~$X$.

(ii) Let $L$ be the extended generator (see, e.g., \cite{RevuzYor}, Definition~7.1.8) of the Markov process $X$.
If the function $V$ belongs to the domain of $L$ and
\[
L V \leq-\varphi\circ V +K,
\]
where $K>0$ and $\varphi\dvtx\mathbb{R}_+\to\mathbb{R}_+$ is a concave
differentiable function increasing to infinity with $\varphi(0)=0$,
then condition~(1) of Theorem~\ref{ThContinuousTime} holds.
\end{Remark}

The proof of this theorem is given in Section~\ref{sect4}. Let us
describe here the main idea.

\begin{pf*}{Sketch of the proof of Theorem~\ref{ThContinuousTime}}
Combining the technique from \cite{DFG,FR,NT}, we find a function
$W\dvtx E\to[0;\infty)$ such that
\[
P^{t_0} W(x)\leq W(x)-\varphi \bigl(K_1W(x) \bigr)
+K_2, \qquad x\in E
\]
for some positive $K_1$, $K_2$. Therefore, by Theorem
\ref{ThDiscreteTime}, the skeleton chain
$(X_{nt_0})_{n\in\mathbb{Z}_+}$ has a unique invariant measure. It is
possible to prove that this measure is also invariant for the Markov
process $X$, and inequality (\ref{CTrate}) holds.
\end{pf*}

Thus Theorems~\ref{ThDiscreteTime} and~\ref{ThContinuousTime} suggest a
new method for proving results concerning subgeometrical convergence.
Namely, one needs to find a suitable contracting metric $d$ and a
suitable Lyapunov function $V$ with $d$-small sublevel sets, such that
the conditions of the theorems hold. It extends the ability of the
existing methods by allowing to choose the metric $d$ (which might be
different from the discrete metric).

\section{Examples and applications}\label{S3}

Let us give some applications of the results of the previous section.
The focus here is on stochastic delay equations; however, it is
possible to apply the results of this kind to study convergence in the
Wasserstein metric for other classes of Markov processes; see, for
example, \cite{HMS}, Section~5.3, for estimates of convergence rates of
stochastic partial differential equations.

We first recall some terminology from \cite{MT}. A Markov chain
$X=(X_n)_{n\in\mathbb{Z}_+}$ is said to be \textit{$\psi$-irreducible}
if there exists a nontrivial measure $\psi$ on $\mathcal{B}(E)$ such
that for any $x\in E$ and any set $A\in\mathcal{B}(E)$ with
$\psi(A)>0$, one has $\P_x(T_A<\infty)>0$, where $T_A$ is the first
return time to the set $A$, that is, $T_{A}:=\inf\{n\geq1\dvtx X_{n}\in
A\}$.

A set $H\in\mathcal{B}(E)$ is called \textit{absorbing} if $P(x,H)=1$
for all $x\in H$, and \textit{Harris} if there exists a measure $\psi$
on $\mathcal{B}(E)$ with $\psi(H)>0$ such that for any $x\in H$ and any
set $A\in\mathcal{B}(E)$ with $\psi(A)>0$ one has $\P_x(T_A<\infty)=1$.

An invariant measure $\pi$ is called \textit{singular} if for any $x\in
E$ there exists an absorbing set $S_x$ such that $x\in S_x$ and
$\pi(S_x)=0$. In other words, the Markov chain, whatever the starting
point is, will remain in the set of $\pi$-measure 0.

\subsection{Autoregressive model}\label{Sect31} Consider the
following peculiar AR(1) process, which belongs to
case (b).

%
\begin{Example}\label{ExAR1}
Let $X=(X_n)_{n\in\mathbb{Z}_+}$ be an autoregressive process
satisfying the following equation:
\[
X_{n+1}=\tfrac1{10}X_n+\varepsilon_{n+1},\qquad n
\in \mathbb{Z}_+,
\]
where $\varepsilon_1, \varepsilon_2, \ldots$ are i.i.d. random
variables uniformly distributed on the set $\{0,\frac1{10},\ldots,
\frac9{10}\}$ and $X_0\in[0;1)$. In other words, to get $X_{n+1}$
from $X_n$ one needs to take the decimal
notation of $X_n$ (which starts with 0 followed by the decimal point)
and insert a random digit immediately
after the decimal point. Other digits in the decimal notation of $X_n$
are shifted right by one position.

Clearly, $X$ is a Markov process with state space $(E,\mathcal
{E})=([0;1),\mathcal{B}([0;1)))$. Let $d$ be the Euclidean metric on
this space [i.e., $d(x,y)=|x-y|$, $x,y\in E$]. One can easily prove
that the process $X$ has a unique invariant measure $\pi$, which is
uniformly distributed on the interval $[0;1)$. Moreover, the sequence
$\{X_n\}$ weakly converges to $\pi$ as $n\to\infty$.
\end{Example}

This autoregression has a number of very interesting and unusual
features. First, it has a reconstruction
property. Namely, if we have just one observation of $X_n$, where the
integer $n$ can be arbitrarily large, then
it is possible to find an initial value $X_0$ with probability $1$ by
the following simple formula: $X_0=\{10^n
X_n\}$, where $\{b \}$ denotes the fractional part of a real $b$. In
other words, one just needs to shift right
the decimal point by $n$ positions and drop all the digits which will
be on the left of the decimal point.

Therefore for $x$, $y\in E$, $x\neq y$, the probability measures
$P(x,\cdot)$ and $P(y,\cdot)$ are singular. Hence the process $X$ has
no nontrivial small sets. On the other hand, the whole state space $E$
is $d$-small. Indeed, it is easily seen that
$W_d(P(x,\cdot),P(y,\cdot))\leq|x-y|/10$, for any $x,y\in E$.

Observe also that the process $X$ is not $\psi$-irreducible, and,
furthermore, it has uncountably many pairwise disjoint absorbing sets.
Indeed, it is sufficient to note that for any $x\in E$ the set
$S_x:=\{y\in E\mid\exists m,n\in\mathbb{Z_+}\dvtx\{10^m y\}=\{10^n
x\}\}$ is absorbing, countable and for $x,y\in E$ either $S_x=S_y$ or
$S_x\cap S_y=\varnothing$. By the same argument, the chain $X$ has no
Harris sets. Since $\pi(S_x)=0$, we see that the measure $\pi$ is
singular.

Finally, let us point out that for any $x\in E$, the sequence
$P^n(x,\cdot)$ does not converge to $\pi$ in total variation [moreover,
$d_{\mathrm{TV}}(P^n(x,\cdot),\pi)=2$ for any positive integer $n$]. On
the other hand, $P^n(x,\cdot)$ converges exponentially to $\pi$ in the
Wasserstein metric [moreover, $W_d(P^n(x,\cdot),\pi)\leq 10^{-n}$ for
any positive integer $n$].

\subsection{Stochastic delay equations}\label{Sect3}
In this subsection we present our results on convergence of SDDEs in
the Wasserstein metric.

Fix $r>0$, positive integers $n$, $m$, and let $\mathcal{C}=\mathcal
{C}([-r;0],\mathbb{R}^n)$ be the space of continuous functions from
$[-r;0]$ to $\mathbb{R}^n$ equipped with the supremum norm $\|\cdot\|$.
Following \cite{HMS}, introduce the following family of metrics on the
space~$\mathcal{C}$:
\[
d_\beta(x,y)=1\wedge\|x-y\|/\beta,\qquad\beta>0.
\]

Consider the stochastic differential delay equation
%
%
\begin{equation}
\label{SDDE} \cases{d X(t)=f(X_t)\,dt+g(X_t)\,dW(t), &
\quad $t\geq0$, \vspace*{2pt}
\cr
X_0=x,}
\end{equation}
where $f\dvtx\mathcal{C}\to\mathbb{R}^n$, $g\dvtx\mathcal{C}\to
\mathbb{R}^{n\times m}$, $W$ is an $m$-dimensional Brownian motion,
$x\in\mathcal{C}$ is the initial condition and as usual we use the
notation $X_t(s):=X(t+s)$, $-r\leq s\leq0$. It is clear that the
process $X=(X_t)_{t\geq0}$ defined on the state space
$(\mathcal{C},\mathcal{B}(\mathcal{C}))$ is Markov.

Throughout this section we assume that the drift and the diffusion
satisfy the following conditions:

\begin{itemize}
\item the drift satisfies a one-sided Lipschitz condition, and the
diffusion is Lipschitz; that is, there exists $K>0$ such that
for any $x,y\in\mathcal{C}$
%
%
\begin{equation}
\label{onesideLip} 2 \bigl\langle f(x) - f(y), x(0) - y(0) \bigr
\rangle^{+} + \bigl\llvert\hspace*{-1pt} \bigl\llvert\hspace*{-1pt}
\bigl\llvert g(x) - g(y) \bigr\rrvert\hspace*{-1pt} \bigr\rrvert\hspace*{-1pt}
\bigr\rrvert^2 \leq K \|x - y\|^2;
\end{equation}

\item the diffusion is nondegenerate; that is, for any $x\in\mathcal
{C}$ the matrix $g(x)$ admits a right
inverse $g^{-1}(x)$ and
%
%
\begin{equation}
\label{nd} \sup_{x\in \mathcal{C}} \bigl\llvert\hspace*{-1pt} \bigl\llvert
\hspace*{-1pt} \bigl\llvert g^{-1}(x) \bigr\rrvert\hspace*{-1pt} \bigr
\rrvert\hspace*{-1pt} \bigr\rrvert<\infty;
\end{equation}

\item\hypertarget{contbbs} (3.4)\hspace*{35pt} $f$ is continuous and bounded on bounded subsets of $\mathcal{C}$.
\end{itemize}
Here $\langle\cdot,\cdot\rangle$ is the standard scalar product in
$\mathbb{R}^n$; for a real $b$ we write $b^+:=\max(b,0)$, and $\llvert\hspace*{-1pt}
\llvert\hspace*{-1pt} \llvert M \rrvert\hspace*{-1pt} \rrvert\hspace*{-1pt} \rrvert$ denotes the Frobenius
norm of a matrix $M$, that is, $\llvert\hspace*{-1pt} \llvert\hspace*{-1pt} \llvert M \rrvert\hspace*{-1pt} \rrvert\hspace*{-1pt}
\rrvert^2=\sum M^2_{ij}$. As in \cite{VerTVP} we also define
%
\[
\lambda_+=\mathop{\sup_{x\in\mathcal{C}}}_{x(0)\neq0}\biggl\langle g(x) g^{T}(x)\frac{x(0)}{|x(0)|},\frac{x(0)}{|x(0)|}
\biggr\rangle,\qquad \Lambda=\sup_{x\in\mathcal{C}}\frac{\Tr g(x) g^{T}(x)}{n}.
\]

\setcounter{equation}{4}
Conditions (\ref{onesideLip}) and \hyperlink{contbbs}{(3.4)} imply \cite{RS} the
existence and uniqueness of the strong solution of SDDE (\ref{SDDE}).

Now we give a general theorem, which describes convergence rates in the
Wasserstein metric~$W_{d_\beta}$. Theorem~\ref{Th3}(i) is a
generalization of \cite{HMS}, Assumption~5.1.

\begin{Theorem}\label{Th3}
Suppose conditions (\ref{onesideLip})--\hyperlink{contbbs}{(3.4)} hold, and there
exists a Lyapunov function $V\dvtx\mathcal C\to\mathbb{R}_+$ that
satisfies inequality (\ref{CTintegr}). If either
\begin{longlist}[(ii)]
\item[(i)] $\lim_{\|x\|\to\infty}V(x)=\infty$
\end{longlist}
or
\begin{longlist}[(ii)]
\item[(ii)] $V(x)=U(x(0))$, for some function
$U\dvtx\mathbb{R}^n\to\mathbb{R}_+$,
$\lim_{|v|\to\infty}U(v)=\infty$, the diffusion coefficient is
uniformly bounded, and the drift coefficient can be decomposed into
two terms,
%
%
\begin{equation}
\label{dec} f(x) = f_1(x) + f_2 \bigl(x(0) \bigr),\qquad x\in\mathcal{C},
\end{equation}
where the function $f_1$ is bounded;
\end{longlist}
then SDDE (\ref{SDDE}) has a unique invariant measure $\pi$.
Furthermore, for any $\beta>0$, the rate of convergence of $\Law(X_t)$
to $\pi$ in the Wasserstein metric $W_{d_\beta}$ is given
by~(\ref{CTrate}).
\end{Theorem}

\begin{pf}
Fix $\beta>0$. Let us check that the process $X$ and the function $V$
satisfy the conditions of Theorem~\ref{ThContinuousTime}. It follows
from \cite{HMS}, Proposition 5.4, and \cite{Shir}, Lemma 3.7.2, that the
process $X$ is Feller. Since $X$ has continuous paths, we see that $X$
is strongly Markovian. The first condition of the theorem is satisfied
by assumption. The second condition also holds. In case (i) it follows
directly from~\cite{HMS}, Section~5.2, that there exists a
$\gamma\in(0;\beta)$ such that the third and the fourth conditions are
met. In case (ii), arguing as in \cite{HMS}, Proposition 5.3 and Lemma
3.8, one can show that the set $\{x\in\mathcal{C}\dvtx|x(0)|\leq R\}$,
$R\geq0$ is $d_\gamma$-small for some $\gamma\in(0;\beta)$, and the
metric $d_\gamma$ is contracting. Thus, in both cases the conditions of
Theorem~\ref{ThContinuousTime} are satisfied.

Apply Theorem~\ref{ThContinuousTime} to the process $X$. It follows
from this theorem that SDDE (\ref{SDDE}) has a unique invariant measure
$\pi$, and the rate of convergence of $\Law(X_t)$ to $\pi$ in the
metric $W_{d_\gamma}$ is provided in (\ref{CTrate}). To complete the
proof, it remains to note that for any measures $\mu_1,
\mu_2\in\mathcal{P}(E)$ one has $W_{d_\beta}(\mu_1,\mu_2)\leq
W_{d_\gamma}(\mu_1,\mu_2)$.
\end{pf}

Ergodic properties of stochastic differential equations (SDE) were
studied by Veretennikov \cite{Ver97,VerTVP}, Malyshkin \cite{Mal},
Klokov \cite{KV}, Douc, Fort and Guillin \cite{DFG} and many others. It
is known that the Veretennikov--Khasminskii condition on the drift
combined with a certain nondegeneracy condition on the diffusion is
sufficient for the existence and uniqueness of the invariant measure
for the strong solution of an SDE. Moreover, these conditions yield
exponential, subexponential or polynomial (depending on the value of
the constant $\alpha$, see below) convergence toward the invariant
measure in the total variation metric \cite{PV,DFG}. The following
theorem extends these results to SDDE.

\begin{Theorem}\label{ThSDDEcoeff}
Suppose conditions (\ref{onesideLip})--\hyperlink{contbbs}{(3.4)} hold, $\Lambda
<\infty$ and the function $f_1$ in decomposition (\ref{dec}) is
bounded.
\begin{longlist}[(ii)]
\item[(i)] Assume additionally that for some constants $\alpha
\in(0,1]$,
$M>0$, $\varkappa>0$, the generalized Veretennikov--Khasminskii condition
holds, that is,
%
%
\begin{equation}
\label{SDDE2} \bigl\langle f(x),x(0) \bigr\rangle\leq- \varkappa \bigl|x(0) \bigr|^{\alpha},
\qquad x\in\mathcal{C}, \bigl|x(0) \bigr|\geq M.
\end{equation}

Then SDDE (\ref{SDDE}) has a unique invariant measure $\pi$, and
$\Law(X_t)$ converges to $\pi$ in the Wasserstein metric $W_{d_\beta}$
subexponentially (if $0<\alpha<1$) or exponentially (if $\alpha=1$);
that is, for any $\beta>0$ there exists positive constants $C_1$ and
$C_2$ such that
%
%
\begin{equation}
\label{th31} \quad W_{d_\beta} \bigl(P^{t}(x,\cdot),\pi \bigr)
\leq C_1 \exp \bigl\{C_1\|x\|^\alpha-C_2
t^{\alpha/(2-\alpha)} \bigr\}, \qquad x\in\mathcal{C}, t>0.
\end{equation}

\item[(ii)] If (\ref{SDDE2}) holds with $\alpha=0$ and $\varkappa>n\Lambda/2$,
then SDDE (\ref{SDDE}) has a unique invariant measure $\pi$, but
$\Law(X_t)$ converges to $\pi$ in the Wasserstein metric
$W_{d_\beta}$ only polynomially; that is, for any $\beta>0$,
$\varepsilon>0$ there exist $C>0$ such that
\[
W_{d_\beta} \bigl(P^{t}(x,\cdot),\pi \bigr)\leq C \bigl(1+\|x
\|^{2+2
\varkappa_0} \bigr)t^{-\varkappa_0+\varepsilon}, \qquad x\in\mathcal{C}, t>0,
\]
where $\varkappa_0=(\varkappa-n\Lambda/2)\lambda_+^{-1}$.
\end{longlist}
\end{Theorem}

\begin{pf}
The proof is based on the application of Theorem~\ref{Th3}(ii) with a
suitable Lyapunov function $V$. (i) Following \cite{KV}, Section~3 (see
also \cite{DFG}, Proposition~5.2), let
$U\dvtx\mathbb{R}^n\to[0;\infty)$ be a twice continuously
differentiable function such that $U(v)=\exp\{k |v|^\alpha\}$ for
$|v|\geq M_0$. The parameters $M_0\geq M$ and $k\geq0$ will be chosen
later. Take $V(x)=U(x(0))$. By Ito's Lemma, for any $x\in\mathcal{C}$
and $t>0$ one has
\begin{eqnarray*}
\E_x V(X_t)&\leq& V(x)+\alpha k \E_x\int
_0^t \II \bigl( \bigl|X(s) \bigr|\geq M_0
\bigr)V(X_s) \bigl|X(s) \bigr|^{\alpha-2} \bigl\langle X(s),f(X_s)
\bigr\rangle\,ds
\\
&&{}+ \frac12 \alpha k \E_x
\\
&&\hspace*{10pt}{}\times \int_0^t \II \bigl(\bigl\llvert
X(s)\bigr\rrvert \geq M_0 \bigr) V(X_s)\bigl|X(s)\bigr|^{\alpha-2}
\bigl(\lambda_+\alpha k \bigl|X(s)\bigr|^{\alpha}+C_1 \bigr)\,ds
\\
&&{}+C_2t
\\
&\leq& V(x) - C_3\alpha k \E_x\int_0^t
\II \bigl(\bigl|X(s)\bigr|\geq M_0 \bigr)V(X_s)\bigl|X(s)\bigr|^{2\alpha-2}
\,ds+C_2t,
\end{eqnarray*}
where $C_1=\lambda_+(\alpha-2)+ n\Lambda$, $C_2>0$, $C_3=\varkappa-\frac
12\lambda_+\alpha k -\frac12C_1M_0^{-\alpha}$ and in the second
inequality we made use of (\ref{SDDE2}).

Let $\varphi\dvtx\mathbb{R}_+\to\mathbb{R}_+$ be a concave
differentiable function with $\varphi(0)=0$ and $\varphi(t)=t (\ln
t)^{(2\alpha-2)/\alpha}$ for $t\geq e^2$. Take
$k=\frac{\varkappa}{2\lambda_+\alpha}$, and
$M_0=(\frac{C_1}\varkappa)^{1/\alpha}\vee(\frac2k)^{1/\alpha}\vee M$. Then
$U(M_0)\geq e^2$ and
\begin{eqnarray*}
\E_x V(X_t)&\leq& V(x) - C_4
\E_x \int_0^t \II \bigl(\bigl|X(s)\bigr|\geq
M_0 \bigr)V(X_s)\bigl|X(s)\bigr|^{2\alpha-2}
\,ds+C_2t
\\
&=& V(x) - C_5 \E_x\int_0^t
\II \bigl(\bigl|X(s)\bigr|\geq M_0 \bigr)\varphi \bigl(V(X_s)
\bigr)\,ds+C_2t
\\
&\leq& V(x) - C_5\E_x \int_0^t
\varphi \bigl(V(X_s) \bigr)\,ds+C_6t,
\end{eqnarray*}
where $C_4:=\alpha k \varkappa/4$, $C_5:=C_4 k^{2/\alpha-2}$ and $C_6>0$. Thus
the function $V$ satisfies inequality (\ref{CTintegr}). Theorem
\ref{Th3}(ii) now yields the existence and the uniqueness of the
invariant measure $\pi$ and implies estimate (\ref{th31}).

(ii) Now let $U(v)=|v|^k$, where $k>2$. We take $V(x)=U(x(0))$ and
proceed as follows:
\begin{eqnarray*}
\E_x V(X_t)&\leq& V(x)+ \frac12 k \E_x\int
_0^t \bigl|X(s)\bigr|^{k-2} \bigl(2 \bigl\langle
X(s),f(X_s) \bigr\rangle+(k-2)\lambda_+ + n\Lambda \bigr)\,ds
\\
&\leq& V(x) - k C_1\E_x\int_0^t
\II \bigl(\bigl|X(s)\bigr|\geq M \bigr) \bigl|X(s)\bigr|^{k-2}\,ds+C_2t,
\end{eqnarray*}
where $C_1= \varkappa-\frac{k-2}2 \lambda_+ - \frac{n\Lambda}{2}$, $C_2>0$. Set
\[
k=2+\frac{2\varkappa-n\Lambda}{\lambda_+}-\varepsilon,
\]
where $\varepsilon>0$. By choosing $\varepsilon>0$ small enough we
can ensure that $k>2$. Take $\varphi(u)=u^{(k-2)/k}$. Then
\[
\E_x V(X_t)\leq V(x) - C_3\E_x
\int_0^t \varphi \bigl(V(X_s)
\bigr)\,ds+C_4t
\]
for some $C_3$, $C_4>0$. Thus the function $V$ satisfies condition
(\ref{CTintegr}), and the statement of the theorem follows now from
Theorem~\ref{Th3}(ii).
\end{pf}

%
\begin{Example}\label{ExSDDE}
Consider the following peculiar SDDE:
\[
dX(t)= f \bigl(X(t) \bigr)\,dt+g \bigl(X(t-1) \bigr)\,dW(t),
\]
where $n=m=1$, the functions $f$ and $g$ satisfies
(\ref{onesideLip})--\hyperlink{contbbs}{(3.4)}, $f$ also satisfies (\ref{SDDE2}),
and $g$ is a strictly increasing bounded positive continuous function.
The strong solution of this SDDE also belongs to case (b). This SDDE
has the reconstruction property \cite{Sch}; that is, if we know $X_t$
for any $t>0$, then we can reconstruct the initial condition $X_0$ with
probability one. Hence, the measures $P^t(x,\cdot)$ and $P^t(y,\cdot)$
are always singular for any $t>0$ and $x\neq y$. It follows from
Theorem~\ref{ThSDDEcoeff} that this SDDE has a unique invariant measure
$\pi$. However, the reconstruction property implies that
$d_{\mathrm{TV}}(P^t(x,\cdot),\pi)$ does not converge to $0$ as
$t\to\infty$, and the measure $\pi$ is singular. On the other hand, if
we replace the total variation metric $d_{\mathrm{TV}}$ by the
Wasserstein metric $W_{d_\beta}$ (these two metrics can be arbitrarily
close to each other for sufficiently small $\beta$), then we see that
$W_{d_\beta}(P^t(x,\cdot),\pi)$ converges to $0$ subexponentially.
\end{Example}
%

\section{Proofs of the main results}\label{sect4}

To prove Theorems~\ref{ThDiscreteTime} and~\ref{ThContinuousTime} we
introduce some notation. Consider a semimetric
$l(x,y):=d(x,y)^{1/p} (1+\beta\varphi(V(x)+V(y)) )^{1/q}$, where
$\beta>0$, $p,q>1$ and $1/p+1/q=1$. These parameters will be chosen
later. We start with two auxiliary lemmas.

\begin{Lemma}\label{LL1}
Assume that a function $V\dvtx E\to[0;\infty)$ satisfies condition (1) of
Theorem~\ref{ThDiscreteTime}. Then for any $n\in\mathbb{Z}_+$
%
%
\begin{equation}
\label{L1first} \sum_{i=0}^{n-1}
P^i ( \varphi\circ V) \leq n K +V.\vadjust{\goodbreak}
\end{equation}

Furthermore, if a measure $\pi$ is invariant for the process $X$, then
$\pi\in\mathcal{P}_{\varphi\circ V}(E)$ and $\pi(\varphi\circ V)\leq
K$.
\end{Lemma}

\begin{pf}
Let us rewrite (\ref{MTLyapunov}) in the following form: $\varphi
\circ V -
K\leq V-PV$. Applying the operator $P^i$, $i\in\mathbb{Z}_+$ to the both
sides of this expression and summing the result over all $0\leq i< n$,
we get
\[
\sum_{i=0}^{n-1} P^i (\varphi
\circ V) - n K \leq V- P^n V,
\]
which proves (\ref{L1first}).

To prove the second part of the lemma we combine the first part of the
lemma with a cut-off argument; see, for example, \cite{Hair06},
Proposition 4.24. Fix $L>0$. Then, for any nonnegative integer $i$, we
have
\begin{eqnarray*}
\int_E \bigl((\varphi\circ V) (x)\wedge L \bigr)
\pi(dx)&=&\int_E P^i \bigl((\varphi\circ V)
\wedge L \bigr) (x) \pi(dx)
\\
&\leq&\int_E \bigl(P^i (\varphi\circ V) (x)
\wedge L \bigr) \pi(dx).
\end{eqnarray*}
Summing the both sides of the above inequality over all $0\leq i< n$,
we derive
\[
\int_E \bigl((\varphi\circ V) (x)\wedge L \bigr) \pi(dx)
\leq\int_E \Biggl( \Biggl(\frac1n\sum
_{i=0}^{n-1} P^i (\varphi\circ V) (x)
\Biggr)\wedge L \Biggr) \pi(dx).
\]
This, combined with (\ref{L1first}), yields
\[
\int_E \bigl((\varphi\circ V) (x)\wedge L \bigr) \pi(dx)
\leq K +\int_E \biggl(\frac
{V(x)}n\wedge L \biggr)
\pi(dx).
\]
Lebesgue's dominated convergence theorem implies that the integral on
the right-hand side of the above
inequality tends to $0$ as $n\to\infty$. Thus
\[
\int_E \bigl((\varphi\circ V) (x)\wedge L \bigr) \pi(dx)
\leq K
\]
and the second part of the lemma follows from Fatou's lemma.
\end{pf}

The following Lemma~\ref{LL2} is due to Petrov.

\begin{Lemma}[(\cite{Petrov})]\label{LL2}
Let $a_0$, $a_1, \ldots$ be a sequence of positive numbers, and assume
that for all $n\in\mathbb{Z}_+$ one has
\[
a_{n+1}\leq a_n \bigl(1-\psi(a_n ) \bigr),
\qquad0\leq a_0\leq1,
\]
where $\psi\dvtx[0;\infty)\to[0;1]$ is a continuous increasing function
with $\psi(0)=0$ and $\psi(x)>0$ for $x>0$. Then
%
%
\begin{equation}
\label{Lconvmain} a_n\leq g^{-1}(n)
\end{equation}
for all $n\in\mathbb{Z}_+$, where
\[
g(x):=\int_x^1 \frac{dt}{t\psi(t)},\qquad0<x
\leq1.
\]
\end{Lemma}

\begin{pf}
We see that the function $g^{-1}$ is well defined. This follows from
the fact that the function $g$ is
nonnegative, unbounded and strictly decreasing. Since $\psi$ is
positive, we have $a_{n+1}\leq a_n$. By the mean
value theorem, there exists $s\in[a_{n+1};a_n]$ such that
\[
g(a_{n+1})-g(a_n)=g'(s) (a_{n+1}-a_n)=-
\frac{a_{n+1}-a_n}{s\psi
(s)}\geq\frac{a_n\psi(a_n)}{s\psi(s)}\geq1.
\]
Hence $g(a_n)\geq n$ and $a_n\leq g^{-1}(n)$.
\end{pf}

The next key lemma gives the estimate of the contraction rate in one step.

\begin{Lemma}\label{LL3}
Assume that the conditions of Theorem~\ref{ThDiscreteTime} hold. Then
there exist $\beta=\beta(p,q)$ and positive $c_1(p,q), c_2(p,q),
c_3(p,q)$ such that for any $\mu,\nu\in\mathcal{P}_{\varphi\circ V}(E)$
one has
\[
W_l(P\mu,P\nu)\leq \bigl(1-c_1\wedge c_2
\varphi' \bigl(\varphi^{-1} \bigl(c W_l( \mu,
\nu)^{-p} \bigr) \bigr) \bigr) W_l(\mu,\nu),
\]
where $c:=c_3(\mu(\varphi\circ V)+\nu(\varphi\circ V))^{p}$ and the
semimetric $l$ was introduced at the beginning of
this section.
\end{Lemma}

Here, as usual, $a\wedge b = \min(a, b)$ and $a\vee b = \max(a, b)$ for
real $a$, $b$. To simplify the formulas, we will drop a pair of
parentheses and write $1-a\wedge b$ for $1-(a\wedge b)$.

\begin{pf*}{Proof of Lemma~\ref{LL3}}
We start as in the proof of \cite{HMS}, Theorem 4.8, by observing that
since $W_l$ is convex, the Jensen inequality implies
%
%
\begin{equation}
\label{MTconvex} W_l(P\mu,P\nu)\leq\int_{E\times E}
W_l \bigl(P(x,\cdot),P(y,\cdot) \bigr) \alpha(dx,dy)
\end{equation}
for any $\mu,\nu\in\mathcal{P}_{\varphi\circ V}(E)$ and any
$\alpha
\in\mathcal{C}(\mu,\nu)$. Applying the
Cauchy--Schwarz inequality and the Jensen inequality for concave
functions, we find that
%
%
\begin{eqnarray}\label{MTstep1}
&& W_l \bigl(P(x,\cdot),P(y,\cdot) \bigr)\nonumber
\\
&&\qquad = \inf_\lambda\int_{E\times E}l(u,v)\lambda(du,dv)\nonumber
\\
&&\qquad \leq\inf_\lambda \biggl( \int_{E\times E}d(u,v) \lambda(du,dv) \biggr)^{1/p}
\\
&&\quad\qquad{}\times \biggl(1+\beta\int_{E\times E} \varphi \bigl(V(u)+V(v) \bigr) \lambda(du,dv)
\biggr)^{1/q}\nonumber
\\
&&\qquad \leq W_d \bigl(P(x,\cdot),P(y,\cdot) \bigr)^{1/p}
\bigl(1+ \beta\varphi \bigl(PV(x)+PV(y) \bigr) \bigr)^{1/q},\nonumber
\end{eqnarray}
where the infimum is taken over all measures $\lambda\in\mathcal
{C}(P(x,\cdot),P(y,\cdot))$.

To estimate the right-hand side of the last inequality we consider
three different cases. Note once again that contrary to the proof of
\cite{HMS}, Theorem 4.8, it is impossible here to obtain a nontrivial
upper uniform bound for $W_l(P(x,\cdot),P(y,\cdot))/l(x,y)$.

Fix a large $M>R$.

\textit{Case} 1. $V(x)+V(y)\leq R$. In this case we proceed similar to
\cite{Hair,HMS}. Using~(\ref{MTstep1}) and conditions (1) and (4) of the
theorem, we obtain
\[
W_l \bigl(P(x,\cdot),P(y,\cdot) \bigr)\leq(1-\rho)^{1/p}
\,d(x,y)^{1/p} \bigl(1+\beta\varphi(2K+R) \bigr)^{1/q}.
\]
Setting
\[
\beta=\frac{(1+\rho/(2-2\rho))^{q-1}-1}{\varphi(2K+R)}
\]
we get
\[
W_l \bigl(P(x,\cdot),P(y,\cdot) \bigr)\leq(1-\rho/2)^{1/p}
\, d(x,y)^{1/p}\leq(1-\rho/2p) l(x,y).
\]

\textit{Case} 2. $R<V(x)+V(y)\leq M$. In this case we make use of
(\ref{MTLyapunov}) and the concavity of $\varphi$ to derive
%
%
\begin{eqnarray}\label{MTcase2}
&& \varphi \bigl(PV(x)+PV(y) \bigr)\nonumber
\\
&&\qquad \leq\varphi \bigl(V(x)+V(y)-
\varphi \bigl(V(x) \bigr)-\varphi \bigl(V(y) \bigr)+2K \bigr)
\\
&&\qquad \leq\varphi \bigl(V(x)+V(y)-\varphi \bigl(V(x)+V(y) \bigr)+2K \bigr).\nonumber
\end{eqnarray}
Clearly, if $u\in(R;M]$, then again by the concavity of $\varphi$ we have
\begin{eqnarray*}
\varphi \bigl(u-\varphi(u)+2K \bigr)&\leq&\varphi(u) \biggl(1- \bigl(\varphi
(u)-2K \bigr)\frac{\varphi
'(u)}{\varphi(u)} \biggr)
\\
&\leq&\varphi(u) \bigl(1-\theta\varphi'(M) \bigr),
\end{eqnarray*}
where $\theta:=1-2K/\varphi(R)$. This inequality, combined with
(\ref{MTstep1}), (\ref{MTcase2}) and contraction property
(\ref{MTdcontr}), yields
\begin{eqnarray*}
&& W_l \bigl(P(x,\cdot),P(y,\cdot) \bigr)
\\
&&\qquad \leq d(x,y)^{1/p}
\bigl(1+\beta\varphi \bigl(PV(x)+PV(y) \bigr) \bigr)^{1/q}
\\
&&\qquad\leq d(x,y)^{1/p} \bigl(1+\beta\varphi \bigl(V(x)+V(y) \bigr)
\bigl(1- \theta\varphi'(M) \bigr) \bigr)^{1/q}
\\
&&\qquad \leq l(x,y) \biggl(1-\frac{\theta\beta\varphi(R)}{1+\beta\varphi
(R)}\varphi'(M)\biggr)^{1/q}
\\
&&\qquad \leq  l(x,y) \biggl(1-\frac{\theta\beta\varphi(R)}{q(1+\beta\varphi
(R))}\varphi'(M) \biggr).
\end{eqnarray*}

\textit{Case} 3. $V(x)+V(y)> M$. This is the easiest situation because
in this case we would like to derive a very weak estimate of
$W_l(P(x,\cdot),P(y,\cdot))$. Combining~(\ref{MTdcontr}),
(\ref{MTstep1}) and (\ref{MTcase2}), we get
\begin{eqnarray*}
&& W_l \bigl(P(x,\cdot),P(y,\cdot) \bigr)
\\
&&\qquad \leq d(x,y)^{1/p} \bigl(1+\beta\varphi \bigl(V(x)+V(y)-\varphi
\bigl(V(x)+V(y) \bigr)+2K \bigr) \bigr)^{1/q}
\\
&&\qquad \leq d(x,y)^{1/p} \bigl(1+\beta\varphi \bigl(V(x)+V(y) \bigr)
\bigr)^{1/q}
\\
&&\qquad =l(x,y).
\end{eqnarray*}
Now we return to the main line of the proof. Introduce
\[
c_1=c_1(p,q,R,K):=\frac{\theta\beta\varphi(R)}{q(1+\beta\varphi
(R))},\qquad
c_2=\rho/2p.
\]
Note that the values of $c_1$ and $c_2$ depend neither on the choice of
$M$ nor on measures $\mu$ and $\nu$. We see from (\ref{MTconvex}) and
the above estimates of $W_l(P(x,\cdot),P(y,\cdot))$ that for all $M>R$
one has
%
%
\begin{eqnarray}\label{MTIntest}
&& W_l(P\mu,P\nu)\nonumber
\\
&&\qquad \leq \bigl(1-c_2\wedge c_1 \varphi'(M) \bigr)\int_{E\times E}
l(x,y) \alpha(dx,dy)
\\
&&\quad\qquad{}+ \bigl(c_2 \wedge c_1 \varphi'(M)
\bigr) \int_{\{V(x)+V(y)>M\}} l(x,y) \alpha(dx,dy).\nonumber
\end{eqnarray}
The second integral on the right-hand side of (\ref{MTIntest}) is
estimated using Chebyshev inequality. Namely,
\begin{eqnarray*}
&&\int_{\{V(x)+V(y)>M\}} l(x,y) \alpha(dx,dy)
\\
&&\qquad\leq\int_{\{V(x)+V(y)>M\}} \bigl(1+\beta\varphi \bigl(V(x)+V(y)
\bigr) \bigr)^{1/q} \alpha(dx,dy)
\\
&&\qquad\leq C \int_{\{V(x)+V(y)>M\}} \varphi \bigl(V(x)+V(y)
\bigr)^{1/q} \alpha(dx,dy)
\\
&&\qquad\leq C \varphi(M)^{-1/p}\int_{E\times E} \varphi
\bigl(V(x)+V(y) \bigr) \alpha(dx,dy)
\\
&&\qquad\leq C \varphi(M)^{-1/p} \bigl(\mu(\varphi\circ V)+\nu(\varphi
\circ V) \bigr),
\end{eqnarray*}
where $C=1/K+\beta+1$, and in the second inequality we used the bound
\mbox{$\varphi(M)>K$}. Note that $\mu(\varphi\circ V)$
as well as $\nu(\varphi\circ V)$ are finite because it was assumed that
$\mu,\nu\in\mathcal{P}_{\varphi\circ V}(E)$.

Recall that $\alpha$ is an arbitrary element of $\mathcal{C}(\mu,\nu
)$. Hence we can take the infimum over all
$\alpha\in\mathcal{C}(\mu,\nu)$ in (\ref{MTIntest}) and use the above
inequality to derive
%
%
\begin{eqnarray}\label{MTnearfinal}
W_l(P\mu,P\nu)&\leq& \bigl(1-c_2
\wedge c_1 \varphi'(M) \bigr)W_l(\mu,\nu)
\nonumber\\[-8pt]\\[-8pt]
&&{}+C \bigl(c_2 \wedge c_1 \varphi'(M)
\bigr) \bigl(\mu(\varphi\circ V)+\nu(\varphi\circ V) \bigr)\varphi
(M)^{-1/p}.\nonumber
\end{eqnarray}
Now we can choose $M$ in such a way, that the right-hand side of the
above expression is always smaller than
$W_l(\mu,\nu)$. Namely, it is sufficient to require that
\[
C \bigl(\mu(\varphi\circ V)+\nu(\varphi\circ V) \bigr)\varphi(M)^{-1/p}
\leq W_l(\mu,\nu)/2.
\]
This inequality holds for
\[
M=\varphi^{-1} \bigl(c_3 \bigl(\mu(\varphi\circ V)+\nu(
\varphi\circ V) \bigr)^{p}W_l(\mu,\nu)^{-p}
\bigr),
\]
where $c_3=c_3(p,q,R,K)=2^{p}(1/K+\beta+1)^{p}$. The substitution of
the last expression into (\ref{MTnearfinal}) proves the lemma.
\end{pf*}

\begin{Lemma}\label{LL4}
Assume that the conditions of Theorem~\ref{ThDiscreteTime} are
satisfied. Let $\mu,\nu\in\mathcal{P}_{\varphi\circ V}(E)$ and let
$(n_k)_{k\in\mathbb{Z}_+}$ be an increasing sequence of positive
integers such that for all $k\in\mathbb{Z}_+$
\[
P^{n_k}\mu(\varphi\circ V)+P^{n_k}\nu(\varphi\circ V)\leq C(
\mu,\nu),
\]
where $C(\mu,\nu)\geq1$. Then there exist positive $C_1,C_2$ that do
not depend on $\mu,\nu$ such that for all
$k\in\mathbb{Z}_+$,
%
%
\begin{equation}
\label{L4main} W_l \bigl(P^{n_k}\mu,P^{n_k}\nu
\bigr)\leq C_1 C(\mu,\nu)\frac{1}{ \varphi
(H_\varphi^{-1}(C_2k))^{1/p}}.
\end{equation}
\end{Lemma}

\begin{pf}
We begin by observing that for any measures $\zeta_1,\zeta_2\in
\mathcal{P}_{\varphi\circ V}(E)$ one has
\begin{eqnarray*}
W_l(\zeta_1,\zeta_2)&\leq&\int
_{E\times E} \bigl(1+\beta\varphi \bigl(V(x)+V(y) \bigr)
\bigr)^{1/q} \zeta_1(dx)\zeta_2(dy)
\\
&\leq& \biggl(1+\beta\int_{E\times E}\varphi \bigl(V(x)+V(y) \bigr)
\zeta_1(dx)\zeta_2(dy) \biggr)^{1/q}
\\
&\leq& \bigl(1+\beta\zeta_1(\varphi\circ V)+\beta
\zeta_2( \varphi\circ V) \bigr)^{1/q},
\end{eqnarray*}
where we used the concavity of the function $\varphi$ and the bound
$d\leq
1$. Hence,
%
%
\begin{eqnarray}
\label{L4a0} W_l \bigl(P^{n_0}\mu,P^{n_0}\nu
\bigr)& \leq &\bigl(1+\beta P^{n_0}\mu(\varphi\circ V) + \beta
P^{n_0} \nu(\varphi\circ V) \bigr)^{1/q}
\nonumber\\[-8pt]\\[-8pt]
&\leq& \bigl(1+\beta C(\mu,\nu) \bigr)^{1/q}\leq(1+\beta)C(\mu,\nu).\nonumber
\end{eqnarray}
Introduce $c_0:=1+\beta$ and denote
\[
a_n:=\frac{W_l(P^n\mu,P^n\nu)}{c_0C(\mu,\nu)},\qquad n\in\mathbb{Z}_+.
\]
It follows from Lemma~\ref{LL3} that $0\leq a_{n+1}\leq a_n$ for all
$n\in\mathbb{Z}_+$. Besides, by definition and (\ref{L4a0}) we have
$a_{n_0}\leq1$. The function $\varphi'$ is decreasing, therefore using
Lemma~\ref{LL3}, we derive
\begin{eqnarray*}
a_{n_{k+1}} &\leq& a_{n_k+1}
\\
&\leq& \bigl(1 - c_1
\\
&&\hspace*{5pt}{}\wedge c_2
\varphi' \bigl(\varphi^{-1} \bigl(c_3
c_0^{-p} \bigl(P^{n_k}\mu(\varphi\circ V) +
P^{n_k}\nu(\varphi\circ V) \bigr)^{p} C(\mu,
\nu)^{-p}a_{n_k}^{-p} \bigr) \bigr) \bigr)
a_{n_k}
\\
&\leq& \bigl(1- c_1\wedge c_2
\varphi' \bigl(\varphi^{-1} \bigl(c_4
a_{n_k}^{-p} \bigr) \bigr) \bigr) a_{n_k},
\end{eqnarray*}
where $c_4=c_0^{-p}c_3$. Since $a_{n_0}\leq1$, it is possible to apply
Lemma~\ref{LL2} to the sequence $(a_{n_k})_{k\in\mathbb{Z}_+}$. It
follows from (\ref{Lconvmain}) that $a_{n_k}\leq g^{-1}(k)$, where
\begin{eqnarray*}
g(x)&=&\int_x^1\frac{dt}{c_1t\wedge c_2 t\varphi'(\varphi^{-1}(c_4 t^{-p}
))}=c_5
\int_x^{c_6}\frac{dt}{t\varphi'(\varphi^{-1}(c_4 t^{-p}
))}+c_7
\\
&=&c_8\int_{c_9}^{\varphi^{-1}(c_4x^{-p})}
\frac{du}{\varphi(u)}+c_7=c_8 H_\varphi \bigl(
\varphi^{-1} \bigl(c_4x^{-p} \bigr)
\bigr)+c_{10}
\end{eqnarray*}
and $c_5,c_6,\ldots $ are some positive constants. Note that to obtain the
third identity, we made the change of
variables $u=\varphi^{-1}(c_4 t^{-p} )$. Thus we finally get
$a_{n_k}\leq
c_{11}\varphi(H_\varphi^{-1}(c_{12}k))^{-1/p}$
and hence
\[
W_l \bigl(P^{n_k}\mu,P^{n_k}\nu \bigr)\leq
c_{13} C(\mu,\nu) \varphi \bigl(H_\varphi^{-1}(c_{12}k)
\bigr)^{-1/p}.
\]
This completes the proof of the lemma.
\end{pf}

\begin{Lemma}\label{LL5}
Under the conditions of Theorem~\ref{ThDiscreteTime}, the process $X$
has a unique stationary measure~$\pi$.
\end{Lemma}

As was pointed out by the referee, if we additionally assumed that the
sublevel sets of $V$ are compact, and the process $X$ is Feller, then
the proof of the lemma would be trivial. Indeed, in this case the
statement of the lemma would follow directly from the
Krylov--Bogoliubov theorem; see \cite{Hair}, page 20. However, we do
not make this assumption because we would like to apply Theorem
\ref{ThDiscreteTime} to Markov processes with a nonlocally compact
state space and in particular, to strong solutions of stochastic delay
equations defined on $\mathcal{C}([-r;0],\mathbb{R}^n)$; see
Section~\ref{Sect3}.

\begin{pf*}{Proof of Lemma~\ref{LL5}}
First let us prove the existence of a stationary measure. Fix $x\in E$.
Let us verify that the sequence of
measures $(P^n \delta_x)_{n\in\mathbb{Z}_+}$ has a~Cauchy
subsequence. For $n<m\in\mathbb{Z}_+$, define
\begin{eqnarray*}
A(n,m)&:=&\# \bigl\{i\in[n;m)\dvtx P^i(\varphi\circ V) (x)\leq4K+4V(x)+1 \bigr\},
\\
B(n,m)&:=&\# \bigl\{i\in[0;n)\dvtx \bigl(P^i(\varphi\circ V) (x)\vee
P^{m-n+i}(\varphi\circ V) (x) \bigr)
\\
&&\hspace*{134pt}\leq4K+4 V(x)+1 \bigr\}.
\end{eqnarray*}
Here the symbol $\#$ denotes the cardinality of a finite set. It
follows from the above definitions that for
$n<m$,
%
%
\begin{equation}
\label{BandA} B(n,m)\geq A(0,n)+A(m-n,n)-n.
\end{equation}

Introduce the following sequence. Let $r_{-1}=-1$ and for $k\in\mathbb{Z}_+$,
\[
r_k:=\inf \bigl\{s>r_{k-1}\dvtx \bigl(P^s(
\varphi \circ V) (x)\vee P^{m-n+s}(\varphi\circ V) (x) \bigr)\leq4K+4 V(x)
\bigr\}.
\]
We see that $r_{B(n,m)-1}<n$. We apply Lemma~\ref{LL4} to the sequence
$(r_k)_{k\in\mathbb{Z}_+}$, the measures $\delta_x$ and
$P^{m-n}\delta_x$ and take $C(\delta_x,P^{m-n}\delta_x)=4K+4 V(x)+1$.
Then, by~(\ref{L4main}),
%
%
\begin{eqnarray}
\label{L5Vazhnoe}
&& W_l \bigl(P^n\delta_x,P^m \delta_x \bigr)\nonumber
\\
&&\qquad\leq W_l \bigl(P^{r_{B(n,m)-1}}\delta
_x,P^{r_{B(n,m)-1}} \bigl(P^{m-n}\delta_x \bigr)\bigr)
\\
&&\qquad \leq C_1 \bigl(4K+4 V(x)+1 \bigr) \varphi \bigl(H_\varphi
^{-1} \bigl(C_2B(n,m)-C_2 \bigr)
\bigr)^{-1/p},\nonumber
\end{eqnarray}
where we used Lemma~\ref{LL3} to obtain the first inequality. Recall
that the constants $C_1,C_2$ are independent of $n,m$.

It follows from (\ref{L1first}) that for any fixed $n$ there exists an
arbitrarily large $m$ such that $A(mn,(m+1)n)\geq3n/4$. Since
$A(0,n)\geq3n/4$, inequality (\ref{BandA}) implies that for any fixed
$n$ there exists an arbitrarily large $m$ such that $B(n,m)\geq n/2$. It
is clear that for all such $m$, one has
\begin{eqnarray*}
&& W_l \bigl(P^n\delta_x,P^m\delta_x \bigr)
\\
&&\qquad \leq C_1 \bigl(4K+4 V(x)+1 \bigr) \varphi
\bigl(H_\varphi^{-1}(C_2n/2-C_2)
\bigr)^{-1/p}=:\Psi(n).
\end{eqnarray*}
It is evident that $\Psi(n)\to0$, as $n\to\infty$.

Now we can construct the desired Cauchy subsequence. We set $n_0=0$,
and for~$k\in\mathbb{Z}_+$,
\[
n_{k+1}:=\inf \bigl\{m>n_{k}\dvtx B(n_k,m)\geq
n_k/2\mbox{ and } \Psi(m)\leq e^{-(k+1)} \bigr\}.
\]
By the above arguments, we see that the sequence $(n_{k})_{k\in\mathbb
{Z}_+}$ is well defined, $B(n_k,n_{k+1})\geq n_k/2$, and $\Psi
(n_{k})\leq
e^{-k}$. Now we claim that the sequence $(P^{n_k}
\delta_x)_{k\in\mathbb{Z}_+}$ is a Cauchy sequence in the space
$(\mathcal{P}(E),W_d)$. Indeed, using (\ref{L5Vazhnoe}) and the
definition of $n_k$ we derive
\begin{eqnarray*}
W_d \bigl(P^{n_k}\delta_x,P^{n_{k+m}}
\delta_x \bigr)&\leq&\sum_{i=k}^{k+m-1}
W_d \bigl(P^{n_i}\delta_x,P^{n_{i+1}}
\delta_x \bigr)
\\
&\leq&\sum_{i=k}^{k+m-1} W_l
\bigl(P^{n_i}\delta_x,P^{n_{i+1}}\delta_x
\bigr)
\\
&\leq&\sum_{i=k}^{k+m-1}
\Psi(n_i) \leq\sum_{i=k}^{k+m-1}
e^{-i} \leq2e^{-k}
\end{eqnarray*}
for all integers $k, m$. Since the space $(\mathcal{P}(E),W_d)$ is
complete (see, e.g., \cite{BK}, Theorem 1.1.3), we see that there
exists a measure $\pi\in\mathcal{P}(E)$ such that
$W_d(P^{n_k}\delta_x, \pi)\to0$.

Let us verify that the measure $\pi$ is stationary, that is, let us
check that $P\pi=\pi$. Note that the metric
$W_d$ is contractive. Indeed, for any $\mu,\nu\in\mathcal P(E)$, we have
\begin{eqnarray*}
W_d(P\mu,P\nu)&\leq&\inf_{\lambda\in\mathcal C(\mu, \nu)} \int
_{E\times E} W_d \bigl(P(x,\cdot),P(y,\cdot) \bigr)
\lambda(dx,dy)
\nonumber
\\
&\leq&\inf_{\lambda\in\mathcal C(\mu, \nu)} \int_{E\times E}d(x,y)
\lambda(dx,dy)
\nonumber
\\
&=&W_d(\mu,\nu),
\end{eqnarray*}
where we used the Jensen inequality and condition (\ref{MTdcontr}).

Therefore, for any $k\in\mathbb{Z}_+$, we obtain
%
%
\begin{eqnarray}
\label{L5pi} W_d(P\pi,\pi)&\leq& W_d \bigl(P
\pi,P^{n_k+1}\delta_x \bigr)\nonumber
\\
&&{}+W_d
\bigl(P^{n_k}\delta_x, P^{n_k+1}\delta_x
\bigr)+W_d \bigl(P^{n_k}\delta_x, \pi \bigr)
\\
&\leq& 2 W_d \bigl(\pi,P^{n_k}\delta_x
\bigr)+W_l \bigl(P^{n_k}\delta_x,
P^{n_k+1}\delta_x \bigr).\nonumber
\end{eqnarray}
The first term on the right-hand side of the last expression tends to
$0$, as \mbox{$k\to\infty$.} To estimate the second term, we observe that if
$n$ is a positive integer,\vadjust{\goodbreak} then $A(0,n)\geq3n/4$ and
$A(1,n+1)\geq3n/4-1$. Therefore, inequality (\ref{BandA}) implies
$B(n,n+1)>n/2-1$. This, combined with (\ref{L5Vazhnoe}), yields
\[
W_l \bigl(P^{n_k}\delta_x,P^{n_k+1}
\delta_x \bigr)\leq C_1 \bigl(4K+4 V(x)+1 \bigr)
\frac{1}{
\varphi(H_\varphi^{-1}(C_2 n_k/2-2C_2))^{1/p}}.
\]
Hence $W_l(P^{n_k}\delta_x,P^{n_k+1}\delta_x)\to0$ as $k\to\infty$,
and we conclude from (\ref{L5pi}) that $W_d(P\pi,\pi)=0$, which implies
the stationarity of the measure $\pi$.

To complete the proof of the lemma it remains to prove the uniqueness
of stationary measure. Suppose that, on the contrary, the process $X$
has two stationary measures $\pi_1$ and $\pi_2$ and $\pi_1\neq\pi_2$.
By Lemma~\ref{LL1}, $\pi_1, \pi_2\in\mathcal{P}_{\varphi\circ
V}(E)$ and
hence $0<W_l(\pi_1,\pi_2)<\infty$. We make use of stationarity of the
measures and Lemma~\ref{LL3} to obtain
\[
W_l(\pi_1,\pi_2)=W_l(P
\pi_1,P\pi_2)<W_l(\pi_1,
\pi_2).
\]
This contradiction proves the lemma.\vadjust{\goodbreak}
\end{pf*}

\begin{pf*}{Proof of Theorem~\ref{ThDiscreteTime}}
It follows from Lemmas~\ref{LL1}~and~\ref{LL5}, that the process
$X$ has a unique stationary measure $\pi\in\mathcal{P}_{\varphi
\circ
V}(E)$ and $\pi(\varphi\circ V)\leq K$. Fix $x\in E$ and consider the
following sequence. Let $n_0=0$ and
\[
n_{k+1}:=\inf \bigl\{m>n_{k}\dvtx P^m (\varphi
\circ V)\leq2K+2V(x)+1 \bigr\},\qquad k\in\mathbb{Z}_+.
\]
We make use of stationarity of $\pi$, the bound $\pi(\varphi\circ
V)\leq
K$ and the definition of~$n_k$ to derive
\[
P^{n_k}\delta_x(\varphi\circ V)+P^{n_k}\pi(
\varphi\circ V)=P^{n_k}\delta_x(\varphi\circ V)+\pi( \varphi
\circ V)\leq3K+2V(x)+1.
\]
Let us apply Lemma~\ref{LL4} to the measures $\delta_x,\pi$, to the
sequence $(n_k)_{k\in\mathbb{Z}_+}$ and take
$C(\delta_x,\pi)=3K+2V(x)+1$. Clearly, $C(\delta_x,\pi)>1$. It follows
from (\ref{L4main}) that
\[
W_l \bigl(P^{n_k}\delta_x,\pi \bigr)\leq
C_1 \bigl(3K+2V(x)+1 \bigr)\frac{1}{ \varphi
(H_\varphi
^{-1}(C_2k))^{1/p}}.
\]
On the other hand, it follows from (\ref{L1first}) that $n_k\leq2k$. To
complete the proof, it remains to take $1/p=1-\varepsilon$ and note
that
\begin{eqnarray*}
W_d \bigl(P^{2k}\delta_x,\pi \bigr)&\leq&
W_l \bigl(P^{2k}\delta_x,\pi
\bigr)=W_l \bigl(P^{2k}\delta_x,P^{2k-n_k}
\pi \bigr)\leq W_l \bigl(P^{n_k}\delta_x,\pi
\bigr)
\\
&\leq& C_1 \bigl(3K+2V(x)+1 \bigr)\frac{1}{ \varphi(H_\varphi
^{-1}(C_2k))^{1/p}}.
\end{eqnarray*}\upqed
\end{pf*}

To switch from discrete time to continuous time and prove Theorem
\ref{ThContinuousTime}, we combine different methods from
\cite{DFG,FR,NT}. First of all for a set $C\in\mathcal{B}(E)$,
introduce the hitting time delayed by $\delta>0$
\[
\tau_C(\delta):=\inf\{t\geq\delta\dvtx X_t\in C\}
\]
and the hitting and return times of the skeleton chain
\begin{eqnarray*}
\sigma_{m, C}&:=&\inf\{n\in\mathbb{Z}_+\dvtx X_{mn}\in C\};
\\
T_{m,C}&:=&\inf\{n\in\mathbb{Z}_+, n\geq1\dvtx X_{mn}\in C\},
\end{eqnarray*}
where $m>0$. Denote for brevity $C_R:=\{x\in E\dvtx V(x)\le R\}$.

\begin{Lemma}\label{LCTL1} If $R>0$ and $\varphi(R)>K$, then under the
conditions of Theorem~\ref{ThContinuousTime}
\[
\E_x \tau_{\{V(x)\leq R\}}(\delta)\leq\frac{\delta\varphi
(R)+V(x)}{\varphi(R)-K}
\]
for all $x\in E$ and $\delta>0$.
\end{Lemma}

\begin{pf}
Fix $L>\delta$. Observe that if $\delta\leq u< \tau_{C_R}(\delta)$, then
by definition $V(X_u)\geq R$. Combining this with (\ref{CTintegr}) we
obtain
\begin{eqnarray*}
\E_x \bigl(\tau_{C_R}(\delta)\wedge L \bigr)&=& \delta+
\E_x \int_\delta^{\tau_{C_R}(\delta)\wedge L}\,du
\\
&\leq&\delta+
\frac{1}{\varphi
(R)}\E_x \int_\delta^{\tau_{C_R}(\delta)\wedge L}
\varphi \bigl(V(X_u) \bigr)\,du
\\
&\leq&\delta+\frac{V(x)+K\E_x (\tau_{C_R}(\delta)\wedge
L)}{\varphi(R)}.
\end{eqnarray*}
Therefore
\[
\E_x (\tau_{C_R(\delta)}\wedge L)\leq\frac{\delta\varphi
(R)+V(x)}{\varphi(R)-K}.
\]
The desired inequality follows now from the Fatou lemma.
\end{pf}

\begin{Lemma}\label{LCTL2}
Let $m>0$. If $R>Km$ and $\varphi(R-Km)>K$, then under the conditions of
Theorem~\ref{ThContinuousTime},
\[
\E_x T_{m, C_{R}}\leq c_1V(x)+c_2,
\qquad x\in E,
\]
where $c_1=c_1(m,R,K)$ and $c_2=c_2(m,R,K)$ are positive functions that
do not depend on $x$.
\end{Lemma}
\begin{pf} The proof of the lemma uses the ideas from the proof of
\cite{FR}, Proposition 22(ii). However, note that we cannot apply this
proposition directly because in contrast to Fort and Roberts, we
assumed neither that the set $\{V(x)\leq R\}$ is petite nor that the
process $X$ is Harris-recurrent with invariant measure.

Introduce $R '<R-Km$ such that $\varphi(R')>K$. The existence of such
$R'$ follows from the conditions of the lemma. Consider the following
sequence of stopping times:
\[
\tau^0:=0,\qquad\tau^1:=\tau_{C_{R'}}(m),\qquad
\tau^n:=\inf \bigl\{t\geq\tau^{n-1}+m\dvtx X_t\in
C_{R'} \bigr\}
\]
and let $M:=\sup_{x\in C_{R'}} \E_x \tau_{C_{R'}}(m)$. By Lemma
\ref{LCTL1},
\[
M\leq\frac{m\varphi(R')+R'}{\varphi(R')-K}.
\]

For $n\in\mathbb{Z}_+$, $n\geq1$ define $Z_n:=\II\{X_{\lceil\tau
^{n}/m\rceil m}\in C_R\}$, where $\lceil b\rceil$ denotes the upper
integer part of a real $b$. By definition,
$Z_{n}\in\mathcal{F}_{\tau^{n+1}}$, where we denote
$\mathcal{F}_t:=\sigma\{X_s, 0\leq s\leq t\}$. We combine\vadjust{\goodbreak} the strong
Markov property, the Chebyshev inequality and (\ref{CTintegr}) to
obtain
%
%
\begin{eqnarray}\label{L7est}
\P(Z_{n}=1 \mid\mathcal{F}_{\tau^n})&=&1-
\E_{X_{\tau^{n}}}\II\{ X_{\lceil\tau^{n}/m\rceil m-\tau^{n}}\notin C_{R}\}
\nonumber
\\
&\geq& 1-\frac{V(X_{\tau^n})+Km }{R}
\\
&\geq&\frac{R-R'-Km}{R}=:\gamma.\nonumber
\end{eqnarray}
It follows from the choice of $R'$ that $\gamma>0$.

Introduce $\eta:=\inf\{n\in\mathbb{Z}_+, n \geq1\dvtx Z_n=1\}$. Using
the strong Markov property, (\ref{L7est}) and following the same lines
as in the proof of \cite{NT}, Lemma 3.1, we get for $n\geq1$ and $x\in
C_{R'}$,
\begin{eqnarray*}
\E_x\tau^n \II(\eta\geq n)&\leq&\E_x
\tau^{n-1}\II(\eta\geq n-1)\E \bigl(\II(Z_{n-1}=0)\mid
\mathcal{F}_{\tau^{n-1}} \bigr)
\\
&&{}+\E_x\II(\eta\geq n-1)\E \bigl(\tau^n-
\tau^{n-1}\mid\mathcal{F}_{\tau^{n-1}} \bigr)
\\
&\leq& (1-\gamma) \E_x\tau^{n-1}\II(\eta\geq n-1) +(1-
\gamma)^{n-1}M.
\end{eqnarray*}
Since $\E_x\tau^0 \II(\eta\geq0)$ is obviously zero, by induction
we establish the following estimate:
\[
\E_x\tau^n \II(\eta\geq n)\leq n M (1-
\gamma)^{n-1},\qquad x\in C_{R'}.
\]
Thus we have
\[
\E_x\tau^\eta\leq\sum_{n=1}^\infty
\E_x\tau^n \II(\eta\geq n)\leq\frac{M}{\gamma^2},\qquad x
\in C_{R'}.
\]
We combine this with Lemma~\ref{LCTL1} to finally obtain
\begin{eqnarray*}
m E_x T_{m, C_{R}}&\leq&\E_x \tau^1 +
\E_x\E_{X_{\tau^1}} \tau^{\eta}+m
\\
&\leq&\frac{m\varphi(R')+V(x)}{\varphi(R)-K}+ \frac{m\varphi
(R')+R'}{\gamma
^2(\varphi(R')-K)}+m
\\
&\leq& c_1 V(x) +c_2
\end{eqnarray*}
for all $x\in E$. This completes the proof of the statement.
\end{pf}

\begin{pf*}{Proof of Theorem~\ref{ThContinuousTime}}
First let us prove that there exist a Lyapunov function $W\dvtx E\to
[0,\infty)$ and positive constants $K_1$, $K_2$ such that
%
%
\begin{equation}
\label{CThelyap} P^{t_0} W(x)\leq W(x)-\varphi \bigl(K_1W(x)
\bigr) +K_2, \qquad x\in E.
\end{equation}
Choose a sufficiently large $R$ (such that the conditions of Lemma
\ref{LCTL2} hold with $m=t_0$), and let
\[
W(x):=\E_x \sum_{k=0}^{\sigma_{t_0, C_R}}
\varphi \bigl(V(X_{kt_0}) \bigr).
\]
It follows from \cite{MT}, Theorem 11.3.5(i) that for $x\in E$
%
%
\begin{equation}
\label{sledMT}\quad  P^{t_0} W(x)= W(x)-\varphi \bigl(V(x) \bigr)+\II(x\in
C_R) \E_x\sum_{k=1}^{T_{t_0, C_R}}
\varphi \bigl(V(X_{kt_0}) \bigr).
\end{equation}
Using an argument similar to that in the proof of \cite{DFG},
Proposition 4.8(i), we obtain for any $L>0$ and $x\in E$,
\begin{eqnarray*}
&& \E_x\sum_{k=1}^{T_{t_0, C_R}\wedge
L}\varphi
\bigl(1+V(X_{kt_0}) \bigr)-\E_x\int_0^{T_{t_0, C_R}\wedge L}
\varphi \bigl(1+V(X_{st_0}) \bigr)\,ds
\\
&&\qquad\leq\frac12\varphi'(1)Kt_0\E_x
(T_{t_0, C_R}\wedge L).
\end{eqnarray*}
Furthermore, using condition (\ref{CTintegr}) and the concavity of the
function $\varphi$, we get for any $x\in E$,
\begin{eqnarray*}
&&\E_x \int_0^{T_{t_0, C_R}\wedge L} \varphi
\bigl(1+V(X_{st_0}) \bigr)\,ds
\\
&&\qquad = \frac{1}{t_0}\E_x \int
_0^{t_0T_{t_0, C_R}\wedge t_0L} \varphi \bigl(1+V(X_{u})
\bigr)\,du
\\
&&\qquad \leq \varphi(1)\E_x (T_{t_0, C_R}\wedge L) +\frac{1}{t_0}
\E_x \int_0^{t_0T_{t_0, C_R}\wedge t_0L} \varphi
\bigl(V(X_{u}) \bigr)\,du
\\
&&\qquad \leq  V(x)/t_0+ \bigl(\varphi(1)+K \bigr)\E_x
(T_{t_0, C_R}\wedge L).
\end{eqnarray*}
Combining this with the previous inequality and using Lemma~\ref{LCTL2}
and Fatou's lemma, we derive for any $x\in E$,
%
%
\begin{eqnarray}\label{itogosumme}
&& \E_x\sum_{k=1}^{T_{t_0, C_R}}
\varphi \bigl(V(X_{kt_0}) \bigr)\nonumber
\\
&&\qquad \leq  V(x)/t_0+ \bigl(
\varphi (1)+K+\varphi'(1)Kt_0 \bigr)\E_x
T_{t_0, C_R}
\nonumber\\[-8pt]\\[-8pt]
&&\qquad \leq  V(x)/t_0+c_3 \bigl(c_1
V(x)+c_2 \bigr)\nonumber
\\
&&\qquad \leq c_4 V(x)+c_5,\nonumber
\end{eqnarray}
where $c_1$ and $c_2$ are defined in Lemma~\ref{LCTL2}, $c_3:=\varphi
(1)+K+\varphi'(1)Kt_0$, $c_4:=1/t_0+c_1c_3$, $c_5=c_2c_3$. Therefore, by
the concavity of $\varphi$,
\[
W(x)\leq\varphi \bigl(V(x) \bigr)+c_4 V(x)+c_5\leq
V(x) \bigl(\varphi'(1)+c_4 \bigr)+\varphi
(1)+c_5.
\]
This bound, together with (\ref{sledMT}) and (\ref{itogosumme}), yields
\[
P^{t_0} W(x)\leq W(x)-\varphi \bigl(c_6 W(x)
\bigr)+c_4R+c_5+c_7
\]
for some positive $c_6$, $c_7$. Hence the function $W$ satisfies
(\ref{CThelyap}).\vadjust{\goodbreak}

Now the statement of Theorem~\ref{ThContinuousTime} follows from the
corresponding statement for discrete time chains. Indeed, the
application of Theorem~\ref{ThDiscreteTime} to the skeleton chain
$(X_{nt_0})_{n\in\mathbb{Z}_+}$ yields the existence of a measure
$\pi$
such that $P^{t_0}\pi=\pi$. Note that for any $0<s<t_0$ the measure
$\pi_s:=P^{s}\pi$ is also invariant for this skeleton chain. Indeed,
$P^{t_0}\pi_s=P^{t_0+s}\pi=P^{s}P^{t_0}\pi=\pi_s$. On the other hand,
Theorem~\ref{ThDiscreteTime} yields uniqueness of the invariant
measure. Thus, $P^{s}\pi=\pi$ and the measure $\pi$ is invariant for
the process $X$. Arguing as in the proof of Lemma~\ref{LL1}, we see
that $\pi(\varphi\circ V)\leq K$.

It follows from Theorem~\ref{ThDiscreteTime} that for any $\varepsilon
>0$ there exist constants $C_1$, $C_2$ such that for all $x\in E$,
$n\in\mathbb{Z}_+$,
\[
W_d \bigl(P^{n t_0}(x,\cdot),\pi \bigr)\leq
\frac{C_1 (1+V(x))}{ \varphi
(H_\varphi
^{-1}(C_2n))^{1-\varepsilon}}.
\]
We combine this with condition (4) of the theorem to conclude that for
any $t>t_0$,
\begin{eqnarray*}
W_d \bigl(P^{t}(x,\cdot),\pi \bigr)&=&W_d
\bigl(P^{t}(x,\cdot),P^{t_0+t-\lfloor
t/t_0\rfloor t_0}\pi \bigr)
\\
&\leq& W_d \bigl(P^{(\lfloor t/t_0\rfloor-1)t_0}(x,\cdot),\pi \bigr)
\\
&\leq&\frac{C_1 (1+V(x))}{ \varphi(H_\varphi
^{-1}(C_3t))^{1-\varepsilon}}
\end{eqnarray*}
for some $C_3>0$. Here $\lfloor b \rfloor$ denotes the lower integer
part of a real $b$. This completes the proof of Theorem
\ref{ThContinuousTime}.
\end{pf*}

\section*{Acknowledgments}
The author is grateful to Professor A.~V. Bulinski and Professor A.~Yu.
Veretennikov for their help and constant
attention to this work. The author also would like to thank Professor
M. Hairer and F.~V. Petrov for useful
discussions and the referee for his valuable comments and suggestions
which helped to improve the quality of the
paper.



\printaddresses

\end{document}